# ON THE BAHADUR REPRESENTATION OF SAMPLE QUANTILES FOR DEPENDENT SEQUENCES

By Wei Biao Wu

*University of Chicago*

We establish the Bahadur representation of sample quantiles for linear and some widely used nonlinear processes. Local fluctuations of empirical processes are discussed. Applications to the trimmed and Winsorized means are given. Our results extend previous ones by establishing sharper bounds under milder conditions and thus provide new insight into the theory of empirical processes for dependent random variables.

**1. Introduction.** Let $(\varepsilon_k)_{k \in \mathbb{Z}}$ be independent and identically distributed (i.i.d.) random variables and let $G$ be a measurable function such that

$$X_n = G(\ldots, \varepsilon_{n-1}, \varepsilon_n) \tag{1}$$

is a well-defined random variable. Clearly $X_n$ represents a huge class of stationary processes. Let $F(x) = \mathbb{P}(X_n \leq x)$ be the marginal distribution function of $X_n$ and let $f$ be its density. For $0 < p < 1$, denote by $\xi_p = \inf\{x : F(x) \geq p\}$ the $p$th quantile of $F$. Given a sample $X_1, \ldots, X_n$, let $\xi_{n,p}$ be the $p$th ($0 < p < 1$) sample quantile and define the empirical distribution function

$$F_n(x) = \frac{1}{n} \sum_{i=1}^{n} \mathbf{1}_{X_i \leq x}.$$

For simplicity we also refer to $\xi_{n,p}$ as the $p$th quantile of $F_n$. In this paper we are interested in finding asymptotic representations of $\xi_{n,p}$. Assuming that $(X_i)_{k \in \mathbb{Z}}$ are i.i.d. and $f(\xi_p) > 0$, Bahadur [1] first established the almost sure result

$$\xi_{n,p} = \xi_p + \frac{p - F_n(\xi_p)}{f(\xi_p)} + O_{\text{a.s.}}[n^{-3/4}(\log n)^{1/2}(\log \log n)^{1/4}], \tag{2}$$

where a sequence of random variables $Z_n$ is said to be $O_{\text{a.s.}}(r_n)$ if $Z_n/r_n$ is almost surely bounded. Refinements of Bahadur's result in the i.i.d. setting were provided by Kiefer in a sequence of papers; see [19, 20, 21]. In particular, Kiefer [19] showed that if $f'$ is bounded in a neighborhood of $\xi_p$ and $f(\xi_p) > 0$, then

$$(3) \quad \limsup_{n \to \infty} \pm \frac{\xi_{n,p} - \xi_p - [p - F_n(\xi_p)]/f(\xi_p)}{n^{-3/4}(\log \log n)^{3/4}} = \frac{2^{5/4} 3^{-3/4} p^{1/2}(1-p)^{1/2}}{f(\xi_p)}$$

almost surely for either choice of sign. Recent contributions can be found in [4, 10].

Extensions of the above results to dependent random variables have been pursued in [26] for $m$-dependent processes, in [27] for strongly mixing processes, in [16] for short-range dependent (SRD) linear processes and in [17] for long-range dependent (LRD) linear processes. The main objective of this paper is to generalize and refine these results for linear and some nonlinear processes.

Sample quantiles are closely related to empirical processes, and the asymptotic theory of empirical processes is then a natural vehicle for studying their limiting behavior. There is a well-developed theory of empirical processes for i.i.d. observations; see, for example, the excellent treatment by Shorack and Wellner [29]. The celebrated Hungarian construction can be used to obtain asymptotic representations of sample quantiles (cf. Chapter 15 in [29]).

Recently there have been many attempts toward a convergence theory of empirical processes for dependent random variables. Such a theory is needed for the related statistical inference. Ho and Hsing [17] and Wu [31] considered the empirical process theory for LRD sequences and obtained asymptotic expansions, while Doukhan and Surgailis [9] considered SRD processes. Instantaneous transforms of Gaussian processes are treated in [7]. Further references on this topic can be found in the recent survey edited by Dehling, Mikosch and Sørensen [6].

For dependent random variables, powerful tools like the Hungarian construction do not exist in general. To obtain comparable results as in the i.i.d. setting, we propose to employ a martingale-based method. The main idea is to approximate sums of stationary processes by martingales. Such approximation schemes act as a bridge which connects stationary processes and martingales. One can then leverage several results from martingale theory, such as martingale central limit theorems, martingale inequalities, the martingale law of the iterated logarithm, and so on, to obtain the desired results. Gordin [13] first applied the martingale approximation method and established a central limit theory for stationary processes; see also [14]. Wu and Woodroofe [37] present some recent developments. Several of its applications on various problems are given in [15, 18, 31, 32, 34].

Historically many limit theorems for dependent random variables have been established under strong mixing conditions. On the other hand, although the martingale approximation-based approach imposes mild and easily verifiable conditions, it nevertheless may allow one to obtain optimal results, in the sense that they may be as sharp as the corresponding ones in the i.i.d. setting.

In this paper, for some SRD linear processes we obtain the following asymptotic representation of sample quantiles:

$$\xi_{n,p} = \xi_p + \frac{p - F_n(\xi_p)}{f(\xi_p)} + O_{\text{a.s.}}[n^{-3/4}(\log \log n)^{3/4}]$$

(cf. Theorem 1), which gives an optimal bound $n^{-3/4}(\log \log n)^{3/4}$ in view of Kiefer's result (3) for i.i.d. random variables. Sample quantiles for LRD processes and some widely used nonlinear processes are also discussed and similar representations are derived. In establishing such asymptotic representations, we also consider the local and global behavior of empirical processes of dependent random variables.

We next introduce the necessary notation. A random variable $\xi$ is said to be in $\mathcal{L}^q$, $q \geq 1$, if $\|\xi\|_q := [\mathbb{E}(|\xi|^q)]^{1/q} < \infty$. Write $\|\cdot\| = \|\cdot\|_2$. Denote the shift process $\mathcal{F}_k = (\ldots, \varepsilon_{k-1}, \varepsilon_k)$ and the projection operator $\mathcal{P}_k \xi = \mathbb{E}(\xi|\mathcal{F}_k) - \mathbb{E}(\xi|\mathcal{F}_{k-1})$, $k \in \mathbb{Z}$. For a sequence of random variables $Z_n$, we say that $Z_n = o_{\text{a.s.}}(r_n)$ if $Z_n/r_n$ converges to 0 almost surely. Write $a_n \sim b_n$ if $\lim_{n \to \infty} a_n/b_n = 1$.

The rest of the paper is structured as follows. Pointwise and uniform Bahadur representations for SRD linear processes are presented in Section 2 and proofs are given in Section 6. LRD processes and nonlinear time series are discussed in Sections 3 and 4, respectively. Applications to the trimmed and Winsorized means are given in Section 5. Section 7 contains proofs and some discussion of results presented in Section 3.

**2. SRD processes.** A causal (one-sided) linear process is defined by $X_k = \sum_{i=0}^{\infty} a_i \varepsilon_{k-i}$, where $\varepsilon_k$ are i.i.d. random variables and $a_k$ are real coefficients such that $X_k$ exists almost surely. The almost sure existence of $X_n$ can be checked by the well-known Kolmogorov three-series theorem (cf. [3]). Let $f_\varepsilon$ and $F_\varepsilon$ be the density and distribution functions of $\varepsilon$, respectively. Recall that $F$ and $F_n$ are the distribution and the empirical distribution functions of $X_n$ and $\xi_p$ is the $p$th quantile of $F$. Without loss of generality let $a_0 = 1$. Define the truncated process by $X_{n,k} = \sum_{j=n-k}^{\infty} a_j \varepsilon_{n-j}$, $k \leq n$, and the *conditional empirical distribution function* by

$$F_n^*(x) = \frac{1}{n} \sum_{i=1}^{n} \mathbb{E}(\mathbf{1}_{X_i \leq x}|\mathcal{F}_{i-1}) = \frac{1}{n} \sum_{i=1}^{n} F_\varepsilon(x - X_{i,i-1}).$$

Throughout this section we assume that

(4) $$\sup_x[f_\varepsilon(x) + |f'_\varepsilon(x)|] < \infty.$$

It is easily seen that (4) implies $\sup_x[f(x) + |f'(x)|] < \infty$ in view of the relation $F(x) = \mathbb{E}[F_\varepsilon(x - \sum_{i=1}^\infty a_i\varepsilon_{k-i})]$ and the Lebesgue dominated convergence theorem. Define the function $\ell_q(n) = (\log\log n)^{1/2}$ if $q > 2$ and $\ell_q(n) = (\log n)^{3/2}(\log\log n)$ if $q = 2$.

THEOREM 1. *Let $X_k = \sum_{i=0}^\infty a_i\varepsilon_{k-i}$ and assume (4), $f(\xi_p) > 0$ and $\mathbb{E}(|\varepsilon_k|^\alpha) < \infty$ for some $\alpha > 0$.*

(a) *If*

(5) $$\sum_{i=n}^\infty |a_i|^{\min(\alpha/q,1)} = O(\log^{-1/q} n)$$

*for some $q > 2$, then* (i) *there exists $C > 0$ such that $\delta_{n,q} = C\ell_q(n)/[f(\xi_p)\sqrt{n}]$ satisfies*

(6) $$F_n(\xi_p + \delta_{n,q}) \geq p \geq F_n(\xi_p - \delta_{n,q}) \qquad \text{almost surely}$$

*and $|\xi_{n,p} - \xi_p| \leq \delta_{n,q}$ almost surely, and* (ii) *the Bahadur representation holds:*

(7) $$\xi_{n,p} = \xi_p + \frac{p - F_n(\xi_p)}{f(\xi_p)} + O_{\text{a.s.}}[n^{-3/4}(\log\log n)^{1/2}\ell_q^{1/2}(n)].$$

(b) *If*

(8) $$\sum_{i=1}^\infty |a_i|^{\min(\alpha/2,1)} < \infty,$$

*then* (i) *and* (ii) *in* (a) *hold for $q = 2$.*

REMARK 1. If $\alpha = 2$, then the process $(X_k)_{k \in \mathbb{Z}}$ has finite variance, and (8) implies that $(X_k)_{k \in \mathbb{Z}}$ is short-range dependent since its covariances are summable.

REMARK 2. If $\alpha > 2$ and there is a $q > 2$ such that (5) holds, then $\sum_{i=n}^\infty |a_i| = O(\log^{-1/\alpha} n)$. The implication is clear if $q < \alpha$. If $q > \alpha$, then $\sum_{i=n}^\infty |a_i|^{\alpha/q} \geq (\sum_{i=n}^\infty |a_i|)^{\alpha/q}$ and we also have $\sum_{i=n}^\infty |a_i| = O(\log^{-1/\alpha} n)$. Therefore, in the case $\alpha > 2$ it suffices to check (5) for the special case $q = \alpha$ instead of verifying it for a whole range of $q > 2$. The condition $\sum_{i=n}^\infty |a_i| = O(\log^{-1/\alpha} n)$ is fairly mild for a linear process being short-range dependent. For example, it is satisfied if $a_n = O(n^{-1}\log^{-1-1/\alpha} n)$.

Assuming that $\mathbb{E}(|\varepsilon_k|^\alpha) < \infty$ for some $\alpha > 0$ and that $|a_n| = O(n^{-\kappa})$ with $\kappa > 1 + 2/\alpha$, Hesse [16] obtained the representation

$$(9) \qquad \xi_{n,p} = \xi_p + \frac{p - F_n(\xi_p)}{f(\xi_p)} + O_{\text{a.s.}}(n^{-3/4+\gamma}),$$

where $\gamma > [\alpha^2(8\kappa - 5) + 2\alpha(10\kappa - 9) - 13]/(4\alpha\kappa - 2\alpha - 2)^2$. In comparison to Hesse's result, our condition (5) only requires $\kappa > \max(1, 2/\alpha)$. If $q > 2$, then the error term (7) is $O_{\text{a.s.}}[n^{-3/4}(\log\log n)^{1/2}\ell_q^{1/2}(n)] = O_{\text{a.s.}}[n^{-3/4}(\log\log n)^{3/4}]$, which gives an optimal bound; see Kiefer's relationship (3). The bound is much better than the one in (9). For example, if $\alpha = 1$ and $\kappa = 3.01$, then Hesse's result (9) gives the error bound $O_{\text{a.s.}}(n^{-0.0031\cdots})$. On the other hand, in Hesse's result $\varepsilon_i$ does not need to have a density.

REMARK 3. It is unclear whether Kiefer's law of the iterated logarithm (3) can be extended to SRD processes. Our result only provides an upper bound. Kiefer's [19] proof involves extremely meticulous analysis and it depends heavily on the i.i.d. assumption. It seems that Kiefer's arguments cannot be directly applied here.

EXAMPLE 1. Suppose that $\varepsilon_i$ is symmetric and its distribution function $F_\varepsilon(x) = 1 - L(x)/x^\alpha$, $x > 0$, where $0 < \alpha \leq 2$ and $L$ is slowly varying at $\infty$. Here a function $L(x)$ is said to be slowly varying at $\infty$ if, for any $\lambda > 0$, $\lim_{x \to \infty} L(\lambda x)/L(x) = 1$. Notice that $\varepsilon_i$ is in the domain of attraction of symmetric $\alpha$-stable distributions.

Assume that $|a_n| = O(n^{-r})$ for some $r > 2/\alpha$. Then for $q \in (2, r\alpha)$, (5) holds. In this case, $\mathbb{E}(|\varepsilon|^\alpha) = 2\alpha \int_0^\infty x^{-1} L(x)\,dx$ may be infinite. However, there exists a pair $(\alpha', q')$ such that $\mathbb{E}(|\varepsilon|^{\alpha'}) < \infty$ and (5) holds for this pair. Actually, one can simply choose $\alpha' < \alpha$ such that $2 < r\alpha'$ and let $q' = (2 + r\alpha')/2$. Then $\sum_{i=n}^\infty |a_i|^{\min(\alpha'/q',1)} = O(n^{1-r\alpha'/q'})$ with $r\alpha'/q' > 1$ and $\mathbb{E}(|\varepsilon|^{\alpha'}) \leq 1 + \int_1^\infty \mathbb{P}(|\varepsilon|^{\alpha'} > u)\,du = 1 + 2\alpha' \int_1^\infty x^{\alpha' - \alpha - 1} L(x)\,dx < \infty$. By Theorem 1 we have the Bahadur representation (7) with the optimal error bound $O_{\text{a.s.}}[n^{-3/4}(\log\log n)^{3/4}]$.

Theorem 1 establishes Bahadur's representation for a single $p \in (0, 1)$. The uniform behavior of $\xi_{n,p} - \xi_p$ over $p \in [p_0, p_1]$, $0 < p_0 < p_1 < 1$, is addressed in Theorem 2. Such results have applications in the study of the trimmed and Winsorized means; see Section 5. Let $\iota_q(n) = (\log n)^{1/q}(\log\log n)^{2/q}$ if $q > 2$ and $\iota_2(n) = (\log n)^{3/2}(\log\log n)$.

THEOREM 2. *Let $X_k = \sum_{i=0}^\infty a_i \varepsilon_{k-i}$. Assume* (4), $\inf_{p_0 \leq p \leq p_1} f(\xi_p) > 0$ *for some $0 < p_0 < p_1 < 1$ and*

$$(10) \qquad \sup_x |f_\varepsilon''(x)| < \infty.$$

In addition, assume that there exist $\alpha > 0$ and $q \geq 2$ such that $\mathbb{E}(|\varepsilon_k|^\alpha) < \infty$ and

$$\sum_{i=1}^{\infty} |a_i|^{\min(\alpha/q,1)} < \infty. \tag{11}$$

Then (i) $\sup_{p_0 \leq p \leq p_1} |\xi_{n,p} - \xi_p| = o_{\text{a.s.}}[\iota_q(n)/\sqrt{n}]$ and (ii) the uniform Bahadur representation holds:

$$\sup_{p_0 \leq p \leq p_1} \left| \xi_{n,p} - \xi_p - \frac{p - F_n(\xi_p)}{f(\xi_p)} \right| = O_{\text{a.s.}}[n^{-3/4}(\iota_q(n) \log n)^{1/2}]. \tag{12}$$

REMARK 4. Generally speaking, (12) cannot be extended to $p_0 = 0$ and/or $p_1 = 1$. The quantity $\xi_{n,p} - \xi_p$ exhibits an erratic behavior as $p \to 0$ or $1$. The extremal theory is beyond the scope of the current paper.

REMARK 5. If $\varepsilon_0$ has finite moments of any order, then under the condition $\sum_{i=1}^{\infty} |a_i| < \infty$, (12) gives the bound $n^{-3/4}(\log n)^{1/2+\eta}$ for any $\eta > 0$.

REMARK 6. The Kiefer–Bahadur theorem asserts that for i.i.d. random variables the left-hand side of (12) has the optimal order $n^{-3/4}(\log n)^{1/2}(\log \log n)^{1/4}$; see Chapter 15 in [29]. Our bound $n^{-3/4}(\iota_q(n) \log n)^{1/2}$ is not sharp. The reason is that we are unable to obtain a law of the iterated logarithm for $\sup_{a \leq x \leq b} |F_n(x) - F(x)|$; see (54) in the proof of Theorem 2 in Section 6.5, where the weaker result $\sup_{a \leq x \leq b} |F_n(x) - F(x)| = o_{\text{a.s.}}[\iota_q(n)/\sqrt{n}]$ is proved. On the other hand, in proving Theorem 1, we are able to establish a law of the iterated logarithm for $F_n(x) - F(x)$ at a *single point* $x$ [cf. Proposition 1 and (i) of Lemma 10], by which the optimal rate $O_{\text{a.s.}}[n^{-3/4}(\log \log n)^{3/4}]$ in (7) can be derived.

**3. LRD processes.** Let the coefficients $a_0 = 1$, $a_n = n^{-\beta}L(n)$, $n \geq 1$, where $1/2 < \beta < 1$ and $L$ is a function slowly varying at infinity; let $X_k = \sum_{i=0}^{\infty} a_i \varepsilon_{k-i}$, where $\varepsilon_k$ are i.i.d. random variables with mean zero and finite variance. By Karamata's theorem (see, e.g., Theorem 0.6 in [25]), the covariances $\gamma(n) = \mathbb{E}(X_0 X_n) \sim C_\beta n^{1-2\beta} L^2(n)$, where $C_\beta = \mathbb{E}(\varepsilon_k^2) \int_0^\infty x^{-\beta}(1+x)^{-\beta} dx$, are not summable and the process is said to be long-range dependent. The asymptotic behavior of LRD processes is quite different from that of SRD ones. We shall apply the empirical process theory developed in [31] and establish Bahadur's representation for long-range dependent processes.

Let $\Psi_n = \sqrt{n} \sum_{k=1}^{n} k^{1/2-2\beta} L^2(k)$ and

$$\sigma_{n,1}^2 = \|n\bar{X}_n\|^2 \sim \frac{C_\beta}{(1-\beta)(3-2\beta)} n^{3-2\beta} L^2(n). \tag{13}$$

By Karamata's theorem, $\Psi_n \sim n^{2-2\beta} L^2(n)/(3/2 - 2\beta)$ if $\beta < 3/4$, $\Psi_n \sim \sqrt{n} L^*(n)$ if $\beta = 3/4$, where $L^*(n) = \sum_{k=1}^{n} L^2(k)/k$ is also a slowly varying function, and $\Psi_n \sim \sqrt{n} \sum_{k=1}^{\infty} k^{1/2-2\beta} L^2(k)$ if $\beta > 3/4$. Let $A_n(\beta) = \Psi_n^2 (\log n)(\log \log n)^2$ if $\beta < 3/4$ and $A_n(\beta) = \Psi_n^2 (\log n)^3 (\log \log n)^2$ if $\beta \geq 3/4$.

THEOREM 3. *Assume* $\inf_{p_0 \leq p \leq p_1} f(\xi_p) > 0$ *for some* $0 < p_0 < p_1 < 1$, $\mathbb{E}(\varepsilon_i^4) < \infty$ *and*

$$\sum_{i=0}^{2} \sup_x |f_\varepsilon^{(i)}(x)| + \int_{\mathbb{R}} |f_\varepsilon'(u)|^2 \, du < \infty. \tag{14}$$

*Let* $b_n = \sigma_{n,1} (\log n)^{1/2} (\log \log n)/n$. *Then*

$$\sup_{p_0 \leq p \leq p_1} \left| \xi_{n,p} - \xi_p - \frac{p - F_n(\xi_p)}{f(\xi_p)} - \frac{\bar{X}_n^2}{2} \frac{f'(\xi_p)}{f(\xi_p)} \right|$$

$$= O_{\text{a.s.}} \left[ b_n^3 + \frac{\sqrt{b_n \log n}}{\sqrt{n}} + \frac{b_n \sqrt{A_n(\beta)}}{n} \right]. \tag{15}$$

The three terms in the $O_{\text{a.s.}}$ bound of (15) have different orders of magnitude for different $\beta$, and correspondingly the term that dominates the bound is different. If $\beta > 7/10$, since $3/2 - 3\beta < -\beta/2 - 1/4$ and $-\beta < -\beta/2 - 1/4$, it is easily seen that $b_n^3 + b_n \sqrt{A_n(\beta)}/n = o[\sqrt{(b_n \log n)/n}]$ in view of $\Psi_n = O[\sqrt{n} L^*(n) + n^{2-2\beta} L^2(n)]$ and $\sqrt{A_n(\beta)} \leq \Psi_n (\log n)^{3/2} (\log \log n)$. Hence the dominant one in the bound of (15) is $O_{\text{a.s.}}[\sqrt{(b_n \log n)/n}]$. On the other hand, if $\beta < 7/10$, then $\sqrt{(b_n \log n)/n} = o[b_n \sqrt{A_n(\beta)}/n]$, $b_n \sqrt{A_n(\beta)}/n \sim C_1 n^{3(1/2-\beta)} L^3(n)(\log n)(\log \log n)^2$ and $b_n^3 \sim C_2 n^{3(1/2-\beta)} L^3(n)(\log n)^{3/2} \times (\log \log n)^3$ for some $0 < C_1, C_2 < \infty$. So $b_n \sqrt{A_n(\beta)}/n = o(b_n^3)$. For the boundary case $\beta = 7/10$, the situation is more subtle since the growth rate of the slowly varying function $L$ is involved. In summary, noting that $\Psi_n = O[\sqrt{n} L^*(n) + n^{2-2\beta} L^2(n)]$, the error bound of (15) is

$$O\{[n^{3(1/2-\beta)} + n^{(1/2-\beta)/2}/n^{1/2} + n^{1/2-\beta}(n^{1/2} + n^{2-2\beta})/n] L_1(n)\}$$

$$= O[n^{\max(-\beta/2-1/4,\ 3/2-3\beta)} L_1(n)] \tag{16}$$

for some slowly varying function $L_1$. This bound is less accurate than the one for the SRD or the i.i.d. counterparts since $\max(-\beta/2 - 1/4,\ 3/2 - 3\beta) > -3/4$ if $\beta < 1$. If $3/4 < \beta < 1$, then the bound is $O_{\text{a.s.}}[n^{-\beta/2-1/4} L_1(n)]$. See Section 7.1 for more discussion on the sharpness of (15) and (16).

In comparison with Bahadur's representations (2) for i.i.d. observations or (7) for short-range dependent processes, (15) has an interesting and different flavor in that it involves the correction term $\frac{1}{2} \bar{X}_n^2 f'(\xi_p)/f(\xi_p)$. More interestingly, this correction term is not needed if $\beta > 5/6$, which includes some LRD processes. Actually, by Lemma 16 in Section 7, $|\bar{X}_n|^2 = o_{\text{a.s.}}(b_n^2)$. Note that $b_n^2 = o(\sqrt{b_n \log n}/\sqrt{n})$ if $\beta > 5/6$. Then the correction term $\frac{1}{2} \bar{X}_n^2 f'(\xi_p)/f(\xi_p)$ can be absorbed into the bound $\sqrt{b_n \log n}/\sqrt{n}$.

If the dependence of the process is strong enough, then we do need the correction $\frac{1}{2}\bar{X}_n^2 f'(\xi_p)/f(\xi_p)$ for a more accurate representation. Specifically, if $\beta \in (1/2, 5/6)$, then $\sqrt{b_n \log n}/\sqrt{n} = o(\sigma_{n,1}^2/n^2)$, $b_n^3 + b_n\sqrt{A_n(\beta)}/n = o(\sigma_{n,1}^2/n^2)$, and as the central limit theorem $n\bar{X}_n/\sigma_{n,1} \Rightarrow N(0,1)$ holds, the correction term has a nonnegligible contribution.

**4. Nonlinear time series.** In the case that $G$ may not have a linear form, we assume that $G$ satisfies the *geometric-moment contraction* (GMC) condition. On a possibly richer probability space, define i.i.d. random variables $\varepsilon'_j, \varepsilon_{i,k}$, $i, j, k \in \mathbb{Z}$, which are identically distributed as $\varepsilon_0$ and are independent of $(\varepsilon_j)_{j \in \mathbb{Z}}$. The process $X_n$ defined in (1) is said to be geometric-moment contracting if there exist $\alpha > 0$, $C = C(\alpha) > 0$ and $0 < r = r(\alpha) < 1$ such that for all $n \geq 0$,

$$(17) \quad \mathbb{E}[|G(\ldots, \varepsilon_{-1}, \varepsilon_0, \varepsilon_1, \ldots, \varepsilon_n) - G(\ldots, \varepsilon'_{-1}, \varepsilon'_0, \varepsilon_1, \ldots, \varepsilon_n)|^\alpha] \leq Cr^n.$$

The process $X'_n := G(\ldots, \varepsilon'_{-1}, \varepsilon'_0, \varepsilon_1, \ldots, \varepsilon_n)$ can be viewed as a coupled version of $X_n$ with the "past" $\mathcal{F}_0 = (\ldots, \varepsilon_{-1}, \varepsilon_0)$ replaced by the i.i.d. copy $\mathcal{F}'_0 = (\ldots, \varepsilon'_{-1}, \varepsilon'_0)$. Here we shall use (17) as our basic assumption for studying the asymptotic behavior of nonlinear time series. Since (17) only imposes the decay rate of the moment of the distance $|X_n - X'_n|$, it is often easily verifiable. In comparison, the classical strong mixing assumptions are typically difficult to check. Recently Hsing and Wu [18] adopted (17) as the underlying assumption and studied the asymptotic behavior of weighted $U$-statistics for nonlinear time series.

Condition (17) is actually very mild as well. Consider the important special class of *iterated random functions* [11], which is recursively defined by

$$(18) \quad X_n = G(X_{n-1}, \varepsilon_n),$$

where $G(\cdot, \cdot)$ is a bivariate measurable function with the Lipschitz constant

$$(19) \quad L_\varepsilon = \sup_{x' \neq x} \frac{|G(x, \varepsilon) - G(x', \varepsilon)|}{|x - x'|} \leq \infty$$

satisfying

$$(20) \quad \mathbb{E}(\log L_\varepsilon) < 0 \quad \text{and} \quad \mathbb{E}[L_\varepsilon^\alpha + |x_0 - G(x_0, \varepsilon)|^\alpha] < \infty$$

for some $\alpha > 0$ and $x_0$. Diaconis and Freedman [8] showed that under (20) the Markov chain (18) admits a unique stationary distribution. Wu and Woodroofe [36] further argued that (20) also implies the geometric-moment contraction (17); see Lemma 3. Some recent improvements are presented in [35]. Under suitable conditions on model parameters, many popular nonlinear time series models such as TAR, RCA and ARCH satisfy (20). Our main result is given next.

THEOREM 4. *Assume* (17), $\sup_x[f(x)+|f'(x)|] < \infty$ *and* $\inf_{p_0 \le p \le p_1} f(\xi_p) > 0$ *for some* $0 < p_0 < p_1 < 1$. *Then*

$$\sup_{p_0 \le p \le p_1} \left| \xi_{n,p} - \xi_p - \frac{p - F_n(\xi_p)}{f(\xi_p)} \right| = O_{\text{a.s.}}(n^{-3/4} \log^{3/2} n). \tag{21}$$

PROOF. For a fixed $\tau > 2$ let $m = \lfloor \omega \log n \rfloor$, where $\omega = \omega_\tau$ is given in Lemma 1 and $\lfloor t \rfloor$ denotes the integer part of $t$; let

$$\tilde{X}_k = G(\ldots, \varepsilon_{k-m-2,k}, \varepsilon_{k-m-1,k}, \varepsilon_{k-m,k}, \varepsilon_{k-m+1}, \varepsilon_{k-m+2}, \ldots, \varepsilon_{k-1}, \varepsilon_k). \tag{22}$$

Our strategy is to replace the "past" $\mathcal{F}_{k-m} = (\ldots, \varepsilon_{k-m-1}, \varepsilon_{k-m})$ in $X_k$ by the i.i.d. copies $(\ldots, \varepsilon_{k-m-2,k}, \varepsilon_{k-m-1,k}, \varepsilon_{k-m,k})$ so that $(\tilde{X}_k)_{k \in \mathbb{Z}}$ is an $m$-dependent process. When $X_n$ is a linear process, Hesse [16] adopted a truncation argument which *forgets* the past $\mathcal{F}_{k-m}$ and approximates $X_k$ by $G_n(\varepsilon_{k-m+1}, \ldots, \varepsilon_{k-1}, \varepsilon_k)$ for some measurable function $G_n$. Clearly the distribution function of $G_n(\varepsilon_{k-m+1}, \ldots, \varepsilon_{k-1}, \varepsilon_k)$ may be different from $F$. Our coupling argument has the advantage that the marginal distribution function of $\tilde{X}_k$ is still $F$. For $j = 1, 2, \ldots, m$, let

$$\tilde{F}_{n,j}(x) = \frac{1}{1 + A_n(j)} \sum_{i=0}^{A_n(j)} \mathbf{1}_{\tilde{X}_{j+im} \le x} \quad \text{and} \quad \tilde{F}_n(x) = \frac{1}{n} \sum_{i=1}^{n} \mathbf{1}_{\tilde{X}_i \le x}, \tag{23}$$

where $A_n(j) = \lfloor n/m \rfloor$ for $1 \le j \le n - m\lfloor n/m \rfloor$ and $A_n(j) = \lfloor n/m \rfloor - 1$ for $1 + n - m\lfloor n/m \rfloor \le j \le m$. Let $A = A_n = n/m$ and $b_n = c\sqrt{\log A}/\sqrt{A}$, where the constant $c$ will be determined later. Let $\tilde{M}_{n,j}(x) = \tilde{F}_{n,j}(x) - F(x)$ and $\tilde{M}_n(x) = \tilde{F}_n(x) - F(x)$. Since $|\tilde{M}_n(x) - \tilde{M}_n(y)| \le \max_{1 \le j \le m} |\tilde{M}_{n,j}(x) - \tilde{M}_{n,j}(y)|$, by Lemma 2 there is a $\delta_\tau > 0$ such that

$$\begin{aligned}
\mathbb{P}&\left[ \sup_{|x-y| \le b_n} |\tilde{M}_n(x) - \tilde{M}_n(y)| > \frac{\delta_\tau (b_n \log A)^{1/2}}{A^{1/2}} \right] \\
&\le \sum_{j=1}^{m} \mathbb{P}\left[ \sup_{|u| \le b_n} |\tilde{M}_{n,j}(x) - \tilde{M}_{n,j}(y)| > \frac{\delta_\tau (b_n \log A)^{1/2}}{A^{1/2}} \right] = mO(A^{-\tau}),
\end{aligned} \tag{24}$$

and similarly $\mathbb{P}[\sup_x |\tilde{M}_n(x)| > \delta_\tau \sqrt{\log A}/\sqrt{A}] = mO(A^{-\tau})$. Since $\tau > 2$, $mA^{-\tau} = O[n^{-\tau}(\log n)^{\tau+1}]$ is summable over $n$. By Lemma 1 and the Borel–Cantelli lemma, we have

$$\begin{aligned}
\sup_{|x-y| \le b_n} &|[F_n(x) - F(x)] - [F_n(y) - F(y)]| \\
&\le \sup_{|x-y| \le b_n} |\tilde{M}_n(x) - \tilde{M}_n(y)| + \frac{2C_\tau \log n}{n} \\
&= \frac{\delta_\tau (b_n \log A)^{1/2}}{A^{1/2}} + \frac{2C_\tau \log n}{n}
\end{aligned} \tag{25}$$

and $\sup_x |F_n(x) - F(x)| \leq \delta_\tau \sqrt{\log A}/\sqrt{A} + C_\tau(n^{-1} \log n)$ almost surely. Now in $b_n = c\sqrt{\log A}/\sqrt{A}$ we choose $c = (2 + \delta_\tau)/[\inf_{p_0 \leq p \leq p_1} f(\xi_p)]$. Then we have

$$\inf_{p_0 \leq p \leq p_1} [F_n(\xi_p + b_n) - p]$$

$$\geq \inf_{p_0 \leq p \leq p_1} [F(\xi_p + b_n) - p] - \sup_{p_0 \leq p \leq p_1} |F_n(\xi_p) - p|$$

$$\quad - \sup_{|x-y| \leq b_n} |[F_n(x) - F(x)] - [F_n(y) - F(y)]|$$

$$\geq b_n \inf_{p_0 \leq p \leq p_1} f(\xi_p) + O(b_n^2) - [\delta_\tau \sqrt{\log A}/\sqrt{A} + C_\tau(n^{-1} \log n)]$$

$$\quad - [\delta_\tau \sqrt{b_n \log A}/\sqrt{A} + 2C_\tau(n^{-1} \log n)]$$

$$> \sqrt{\log A}/\sqrt{A}$$

almost surely. Similarly $\sup_{p_0 \leq p \leq p_1} [F_n(\xi_p - b_n) - p] < 0$ almost surely. Hence for $\Delta_{n,p} = \xi_{n,p} - \xi_p$, $\sup_{p_0 \leq p \leq p_1} |\Delta_{n,p}| \leq b_n$ almost surely since $F_n$ is nondecreasing. Since $|F_n(\xi_{n,p}) - p| \leq 1/n$, by (25)

$$\sup_{p_0 \leq p \leq p_1} |[F_n(\xi_{n,p}) - F(\xi_{n,p})] - [F_n(\xi_p) - F(\xi_p)]|$$

$$= \sup_{p_0 \leq p \leq p_1} |[p - F(\xi_{n,p})] - [F_n(\xi_p) - F(\xi_p)]| + O(1/n)$$

$$= O_{\text{a.s.}} \left[ \frac{\delta_\tau (b_n \log A)^{1/2}}{A^{1/2}} + \frac{2C_\tau \log n}{n} \right] + O(1/n)$$

$$= O_{\text{a.s.}}(n^{-3/4} \log^{3/2} n),$$

which entails (21) in view of $\inf_{p_0 \leq p \leq p_1} f(\xi_p) > 0$ and, by Taylor's expansion, $F(\xi_{n,p}) - F(\xi_p) = \Delta_{n,p} f(\xi_p) + O(\Delta_{n,p}^2)$ since $\sup_x |f'(x)| < \infty$. □

LEMMA 1. *Assume* (17) *and* $\sup_x f(x) < \infty$. *Then for any* $\tau > 1$, *there exist* $\omega_\tau, C_\tau > 0$ *such that for* $m = \lfloor \omega_\tau \log n \rfloor$ *we have*

(26) $$\mathbb{P}\left[\sup_x |\tilde{F}_n(x) - F_n(x)| \geq C_\tau n^{-1} \log n\right] = O(n^{-\tau}).$$

PROOF. Let $\rho = r^{1/(2\alpha)}$, $\omega_\tau = -(1 + \alpha^{-1})(\tau + 2)/\log \rho$ and $C_\tau = 1 - (1 + \alpha^{-1})(\tau + 1)/\log \rho$; let $R_n$ be the set $\bigcap_{i=1}^n \{|X_i - \tilde{X}_i| \leq \rho^m\}$ and let $R_n'$ be its complement. Then

$$\mathbb{P}(R_n') \leq n\mathbb{P}(|X_i - \tilde{X}_i| \geq \rho^m) \leq n\rho^{-\alpha m} \mathbb{E}(|X_i - \tilde{X}_i|^\alpha)$$

$$\leq n\rho^{-\alpha m} Cr^m = nC\rho^{\alpha m} = o(n^{-\tau}).$$

Let $K = C_\tau - 1$. By the triangle inequality, to establish (26) it suffices to show that

$$(27) \quad \mathbb{P}\left[\sup_x |\tilde{F}_n(x) - F_n(x)| \mathbf{1}_{R_n} > K n^{-1} \log n\right] = O(n^{-\tau}).$$

Notice that $\sup_x |\tilde{F}_n(x) - F_n(x)| \mathbf{1}_{R_n} \leq \sup_x [F_n(x+\rho^m) - F_n(x-\rho^m)]$. Clearly, the event $\{\sup_x [F_n(x+\rho^m) - F_n(x-\rho^m)] > K n^{-1} \log n\}$ implies that there exist two indices $i$ and $j$ with $j - i \geq \lfloor K \log n \rfloor$ such that both $X_i$ and $X_j$ are in the interval $[x - \rho^m, x + \rho^m]$ for some $x \in \mathbb{R}$. Therefore

$$\mathbb{P}\left[\sup_x [F_n(x+\rho^m) - F_n(x-\rho^m)] > K n^{-1} \log n\right]$$

$$\leq \mathbb{P}\left[\bigcup_{i=1}^{n-\lfloor K \log n \rfloor} \bigcup_{j=i+\lfloor K \log n \rfloor}^{n} \{|X_i - X_j| \leq 2\rho^m\}\right]$$

$$\leq \sum_{i=1}^{n-\lfloor K \log n \rfloor} \sum_{j=i+\lfloor K \log n \rfloor}^{n} \mathbb{P}(|X_i - X_j| \leq 2\rho^m)$$

$$\leq n \sum_{j=\lfloor K \log n \rfloor}^{n} \mathbb{P}(|X_0 - X_j| \leq 2\rho^m).$$

Recall (17) for $X'_j = G(\ldots, \varepsilon'_{-1}, \varepsilon'_0, \varepsilon_1, \ldots, \varepsilon_j)$. Then

$$\mathbb{P}(|X_0 - X_j| \leq 2\rho^m)$$

$$\leq \mathbb{P}(|X_0 - X_j| \leq 2\rho^m, |X_j - X'_j| \leq \rho^j) + \mathbb{P}(|X_j - X'_j| > \rho^j)$$

$$\leq \mathbb{P}(|X_0 - X'_j| \leq 2\rho^m + \rho^j) + \rho^{-\alpha j} C r^j.$$

Observe that $X_0$ and $X'_j$ are i.i.d. and $\mathbb{P}(|X_0 - X'_j| \leq \delta) = \mathbb{E}[\mathbb{P}(|X_0 - X'_j| \leq \delta | X'_j)] \leq 2c\delta$, where $c = \sup_x f(x) < \infty$. Thus

$$\mathbb{P}\left[\sup_x |\tilde{F}_n(x) - F(x)| \mathbf{1}_{R_n} > K n^{-1} \log n\right] \leq n \sum_{j=\lfloor K \log n \rfloor}^{n} [2c(2\rho^m + \rho^j) + \rho^{\alpha j} C]$$

$$= nO(n\rho^m + \rho^{K \log n}) + nO(\rho^{\alpha K \log n}),$$

which ensures (27) by the choice of $K$ and $\omega_\tau$. □

LEMMA 2. Let $(Z_k)_{k \in \mathbb{Z}}$ be i.i.d. random variables with distribution and density functions $F_Z$ and $f_Z$ for which $\sup_z f_Z(z) < \infty$; let $F_{n,Z}(z) = \frac{1}{n} \sum_{i=1}^n \mathbf{1}_{Z_i \leq z}$. Then for all $\tau > 1$ there exists $C_\tau > 0$ such that

$$(28) \quad \mathbb{P}\left[\sup_x |F_{n,Z}(x) - F_Z(x)| > \frac{C_\tau (\log n)^{1/2}}{n^{1/2}}\right] = O(n^{-\tau})$$

*and*

$$\mathbb{P}\left[\sup_{|x-y|\leq b_n} |F_{n,Z}(x) - F_Z(x) - \{F_{n,Z}(y) - F_Z(y)\}| > \frac{C_\tau (b_n \log n)^{1/2}}{n^{1/2}}\right]$$
(29)
$$= O(n^{-\tau}),$$

*where $(b_n)_{n\geq 1}$ is a positive, bounded sequence of real numbers such that $\log n = o(nb_n)$.*

Lemma 2 easily follows from classical results for i.i.d. uniform random variables under quantile transformations; see the Dvoretzky–Kiefer–Wolfowitz inequality and Inequality 14.0.9 in [29]. The lemma is needed in the proof of Theorem 4 and it is a special case of Lemma 7 in Section 6.2. We purposely state Lemma 2 here also for the sake of comparison: the martingale-based method may yield results comparable to those obtained under the i.i.d. assumption.

**5. Trimmed and Winsorized means.** Let $\xi_{n,1/n} \leq \xi_{n,2/n} \leq \cdots \leq \xi_{n,1}$ be the order statistics of $X_1, \ldots, X_n$. Then the trimmed and Winsorized means are of the forms $\sum_{i=\alpha(n)+1}^{\beta(n)} \xi_{n,i/n}/[\beta(n) - \alpha(n)]$ and $n^{-1}[\alpha(n)\xi_{n,\alpha(n)/n} + (n - \beta(n))\xi_{n,\beta(n)/n+1/n} + \sum_{i=\alpha(n)+1}^{\beta(n)} \xi_{n,i/n}]$, respectively, where $\alpha(n) = \lfloor np_0 \rfloor$ and $\beta(n) = \lfloor np_1 \rfloor$.

Stigler [30] studied the asymptotic behavior of trimmed means for i.i.d. random variables. Here we shall apply Theorems 2 and 4 to obtain a central limit theorem for some dependent random variables. SRD linear processes and causal processes satisfying (17) are considered in (i) and (ii) of Theorem 5, respectively. Denote by $N(\mu, \sigma^2)$ a normal distribution with mean $\mu$ and variance $\sigma^2$.

THEOREM 5. (i) *Let $q = 2$ and assume that the conditions of Theorem 2 are satisfied. Then there is a $\sigma < \infty$ such that*

$$\sqrt{n}\left[\frac{\sum_{i=\alpha(n)+1}^{\beta(n)} \xi_{n,i/n}}{\beta(n) - \alpha(n)} - \frac{1}{p_1 - p_0}\int_{p_0}^{p_1} \xi_u \, du\right] \Rightarrow N(0, \sigma^2).$$
(30)

(ii) *Assume that the conditions of Theorem 4 are satisfied. Then the central limit theorem (30) holds.*

PROOF. (i) Since $\xi_{n,u}$ is nondecreasing in $u$, $n\int_{(i-1)/n}^{i/n} \xi_{n,u}\,du \leq \xi_{n,i/n} \leq n\int_{i/n}^{(i+1)/n} \xi_{n,u}\,du$ holds for $1 < i < n-1$. Hence

$$n\int_{\alpha(n)/n}^{\beta(n)/n} \xi_{n,u}\,du \leq \sum_{i=\alpha(n)+1}^{\beta(n)} \xi_{n,i/n} \leq n\int_{[1+\alpha(n)]/n}^{[1+\beta(n)]/n} \xi_{n,u}\,du.$$

It is easily seen that, under the conditions of Theorem 2, (12) also holds over the expanded interval $[p_0 - \tau, p_1 + \tau]$ for some sufficiently small $\tau > 0$. Therefore, we have $\sup_{\alpha(n)/n \leq u \leq [1+\beta(n)]/n} |\xi_{n,u}| = O_{\text{a.s.}}(1)$ and consequently

$$\sum_{i=\alpha(n)+1}^{\beta(n)} \xi_{n,i/n} - n \int_{p_0}^{p_1} \xi_{n,u}\, du = O_{\text{a.s.}}(1). \tag{31}$$

By (12) of Theorem 2,

$$\int_{p_0}^{p_1} \xi_{n,u}\, du - \int_{p_0}^{p_1} \xi_u\, du - \int_{p_0}^{p_1} \frac{u - F_n(\xi_u)}{f(\xi_u)}\, du$$
$$= O_{\text{a.s.}}[n^{-3/4}(\iota_2(n)\log n)^{1/2}]. \tag{32}$$

Lemma 11 in Section 6.4 asserts that $\{\sqrt{n}[F_n(x) - F(x)],\ \xi_{p_0} \leq x \leq \xi_{p_1}\} \Rightarrow \{W(x),\ \xi_{p_0} \leq x \leq \xi_{p_1}\}$ for some centered Gaussian process $W$ in the Skorohod space $D[\xi_{p_0}, \xi_{p_1}]$ [2]. By the continuous mapping theorem, (30) follows from (31) and (32).

(ii) By Theorem 4 in [35], under the conditions (17) and $\sup_x f(x) < \infty$, we also have the functional central limit theorem $\{\sqrt{n}[F_n(x) - F(x)],\ \xi_{p_0} \leq x \leq \xi_{p_1}\} \Rightarrow \{W(x),\ \xi_{p_0} \leq x \leq \xi_{p_1}\}$ for some Gaussian process $W$. So (30) holds in view of the argument in (i). □

REMARK 7. Using the same argument, it is easily seen that for the Winsorized mean $n^{-1}[\alpha(n)\xi_{n,\alpha(n)/n} + (n - \beta(n))\xi_{n,\beta(n)/n+1/n} + \sum_{i=\alpha(n)+1}^{\beta(n)} \xi_{n,i/n}]$, we also have the central limit theorem (30) with the asymptotic mean $(p_1 - p_0)^{-1} \int_{p_0}^{p_1} \xi_u\, du$ replaced by $p_0 \xi_{p_0} + (1-p_1)\xi_{p_1} + \int_{p_0}^{p_1} \xi_u\, du$. Other forms of linear functions of order statistics can be similarly handled.

**6. Proofs of Theorems 1 and 2.** We first introduce our method. Recall $\mathcal{F}_k = (\ldots, \varepsilon_{k-1}, \varepsilon_k)$ and $F_n^*(x) = \sum_{i=1}^n F_\varepsilon(x - X_{i,i-1})/n$. Write $F_n(x) - F(x) = M_n(x) + N_n(x)$, where $M_n(x) = F_n(x) - F_n^*(x)$ and $N_n(x) = F_n^*(x) - F(x)$.

Notice that under (4) the conditional empirical distribution function $F_n^*$ is differentiable with the uniformly bounded derivative $f_n^*(x) = n^{-1}\sum_{i=1}^n f_\varepsilon(x - X_{i,i-1})$ and hence $dN_n(x)/dx = f_n^*(x) - f(x)$ is also uniformly bounded. The differentiability property greatly facilitates the related analysis. In comparison, $F_n$ is a step function and hence discontinuous. On the other hand, $nM_n(x)$ forms a martingale with bounded, stationary and ergodic increments $\mathbf{1}_{X_i \leq x} - \mathbb{E}(\mathbf{1}_{X_i \leq x}|\mathcal{F}_{i-1})$. Therefore, results from martingale theory are applicable.

The martingale part $M_n$ and the differentiable part $N_n$ are treated in Sections 6.2 and 6.3, respectively. Section 6.4 discusses the oscillatory behavior and some asymptotic properties of empirical processes, which are needed for the derivation of Bahadur's representations. Proofs of Theorems 1 and 2 are given in Section 6.5.

6.1. *Some useful results.* The following Proposition 1 is needed in proving Theorems 1 and 2. See [33] for a proof.

PROPOSITION 1. *Let $S_n(g) = \sum_{i=1}^n g(\mathcal{F}_i)$, where $g$ is a measurable function such that $g(\mathcal{F}_0) \in \mathcal{L}^q$ for some $q \geq 2$, $\mathbb{E}[g(\mathcal{F}_0)] = 0$ and*

$$(33) \qquad \Theta_{0,q} := \sum_{i=0}^\infty \|\mathcal{P}_0 g(\mathcal{F}_i)\|_q < \infty.$$

*Let $B_q = 18q^{3/2}(q-1)^{1/2}$ if $q > 2$ and $B_q = 1$ if $q = 2$. Then*

$$(34) \qquad \|S_n(g)\|_q \leq B_q \sqrt{n} \Theta_{0,q}.$$

*Furthermore, if $\Theta_{m,q} := \sum_{i=m}^\infty \|\mathcal{P}_0 g(\mathcal{F}_i)\|_q = O[(\log m)^{-1/q}]$ for some $q > 2$, then*

$$(35) \qquad \limsup_{n\to\infty} \pm \frac{S_n(g)}{\sqrt{2n \log \log n}} = \sigma$$

*almost surely for either choice of sign, where $\sigma = \|\sum_{i=0}^\infty \mathcal{P}_0 g(\mathcal{F}_i)\| < \infty$.*

In order to apply Proposition 1 to $S_n(g) = n[F_n(x) - F(x)]$ or $n[f_n^*(x) - f(x)]$, one needs to estimate $\|\mathcal{P}_0 \mathbf{1}_{X_i \leq x}\|$ or $\|\mathcal{P}_0 f_\varepsilon(x - X_{i,i-1})\|$. The following Lemma 3 provides a simple upper bound if the random variable $\varepsilon_0$ satisfies certain moment conditions. In particular, $\varepsilon_0$ is allowed to have infinite variance.

LEMMA 3. *Let $X_k = \sum_{i=0}^\infty a_i \varepsilon_{k-i}$, where $\varepsilon_k$ are i.i.d. with $\mathbb{E}(|\varepsilon_k|^\alpha) < \infty$ for some $\alpha > 0$. Then under (4), $\|\mathcal{P}_0 g(\mathcal{F}_n)\|_q = O[|a_n|^{\min(\alpha/q,1)}]$ holds for $g(\mathcal{F}_n) = \mathbf{1}_{X_n \leq x}$ and $g(\mathcal{F}_n) = f_\varepsilon(x - X_{n,n-1})$. If additionally (10) is satisfied, then the same bound also holds for $g(\mathcal{F}_n) = f'_\varepsilon(x - X_{n,n-1})$.*

PROOF. Let $(\varepsilon'_i)_{i \in \mathbb{Z}}$ be an i.i.d. copy of $(\varepsilon_i)_{i \in \mathbb{Z}}$ and $X_n^* = X_n - a_n \varepsilon_0 + a_n \varepsilon'_0$; let $G_n$ be the distribution function of $X_n - X_{n,0} = \sum_{j=0}^{n-1} a_j \varepsilon_{n-j}$. Since $c = \sup_x f_\varepsilon(x) < \infty$, it is easily seen that the density $g_n(x) = G'_n(x)$ is also bounded by $c$. Observe that $\mathbb{P}(X_n^* \leq x | \mathcal{F}_0) = \mathbb{P}(X_n \leq x | \mathcal{F}_{-1})$. By Jensen's inequality,

$$\|\mathcal{P}_0 g(\mathcal{F}_n)\|_q \leq \|\mathbb{P}(X_n \leq x | \mathcal{F}_0) - \mathbb{P}(X_n^* \leq x | \mathcal{F}_0)\|_q$$
$$= \|\mathbb{E}[G_n(x - X_{n,0}) - G_n(x - X_{n,0} + a_n \varepsilon_0 - a_n \varepsilon'_0) | \mathcal{F}_0]\|_q$$
$$\leq \|G_n(x - X_{n,0}) - G_n(x - X_{n,0} + a_n \varepsilon_0 - a_n \varepsilon'_0)\|_q$$

$$\leq \|\min(c|a_n\varepsilon_0 - a_n\varepsilon_0'|, 1)\|_q$$
$$\leq [\mathbb{E}(c|a_n\varepsilon_0 - a_n\varepsilon_0'|)^{\min(\alpha,q)}]^{1/q}$$
$$= O[|a_n|^{\min(\alpha/q,1)}].$$

Here the elementary inequality $[\min(|b|,1)]^q \leq |b|^{\min(\alpha,q)}$ is applied. The other cases $g(\mathcal{F}_n) = f_\varepsilon(x - X_{n,n-1})$ and $g(\mathcal{F}_n) = f_\varepsilon'(x - X_{n,n-1})$ can be similarly proved. □

To establish a uniform Bahadur representation for $\xi_{n,p} - \xi_p$ over $p \in [p_0, p_1]$, $0 < p_0 < p_1 < 1$, we need the following version of maximal inequality, which will be used to obtain an almost sure upper bound of $\sup_{\xi_{p_0} \leq x \leq \xi_{p_1}} |F_n(x) - F(x)|$. Similar versions appeared in [2, 22, 24, 28]. For a proof of Lemma 4 see [33].

LEMMA 4. *Let $(Y_{k,\theta}, k \in \mathbb{Z})_{\theta \in \Theta}$ be a class of centered stationary processes in $\mathcal{L}^q$, $q > 1$. Namely, for each $\theta \in \Theta$, $(Y_{k,\theta})_{k \in \mathbb{Z}}$ is a stationary process in $\mathcal{L}^q$ and $\mathbb{E}(Y_{k,\theta}) = 0$. Let $S_{n,\theta} = Y_{1,\theta} + \cdots + Y_{n,\theta}$ and let $d = d(n)$ be an integer such that $2^{d-1} < n \leq 2^d$. Then*

$$(36) \quad \left\{\mathbb{E}^*\left[\max_{k \leq n} \sup_{\theta \in \Theta} |S_{k,\theta}|^q\right]\right\}^{1/q} \leq \sum_{j=0}^{d} 2^{(d-j)/q} \left\{\mathbb{E}^*\left[\sup_{\theta \in \Theta} |S_{2^j,\theta}|^q\right]\right\}^{1/q},$$

*where $\mathbb{E}^*$ is the outer expectation $\mathbb{E}^* Z = \inf\{\mathbb{E}X : X \geq Z, X \text{ is a random variable}\}$.*

6.2. *The martingale part $M_n$.*

LEMMA 5. *Let $(b_n)_{n \geq 1}$ be a positive, bounded sequence of real numbers such that $\log^3 n = o(nb_n)$. Assume $\sup_x f_\varepsilon(x) < \infty$. Then for any $\tau > 1$ there exists a constant $C_\tau > 0$ such that*

$$(37) \quad \mathbb{P}\left\{\sup_{|u| \leq b_{2^k}} \max_{2^{k-1} < n \leq 2^k} n|M_n(x+u) - M_n(x)| > C_\tau \sqrt{2^k b_{2^k} \log k}\right\}$$
$$= O(k^{-\tau}).$$

PROOF. Let $c = \sup_x f_\varepsilon(x) < \infty$. For a given $u > 0$, since $\mathbb{P}(x < X_i \leq x + u | \mathcal{F}_{i-1}) \leq cu$, we have

$$\sum_{i=1}^{n} [\mathbb{E}(\mathbf{1}_{x < X_i \leq x+u} | \mathcal{F}_{i-1}) - \mathbb{E}^2(\mathbf{1}_{x < X_i \leq x+u} | \mathcal{F}_{i-1})] \leq ncu.$$

Here without loss of generality we restrict $u$ to be nonnegative. Let $t_k = \sqrt{2^k b_{2^k} \log k}$. Since $\mathbf{1}_{x < X_i \leq x+u} - \mathbb{E}(\mathbf{1}_{x < X_i \leq x+u} | \mathcal{F}_{i-1})$, $1 \leq i \leq n$, form bounded martingale differences, by Freedman's inequality (cf. Theorem 1.6 in [12]) we get that

$$\mathbb{P}\left\{\max_{2^{k-1}<n\leq 2^k} n|M_n(x+u) - M_n(x)| > Ct_k\right\}$$

(38)
$$\leq 2\exp[-C^2 t_k^2/(2Ct_k + 2\times 2^k cu)]$$

for all $C > 0$. Let $\alpha_k = b_{2^k}/k$, $u_i = i\alpha_k$, $i = 0, 1, \ldots, k-1$, and $v_m = mb_{2^k}/(k2^k)$, $m = 0, 1, \ldots, 2^k - 1$. Since $t_k = o(2^k b_{2^k}/k)$, we have for sufficiently large $k$ that

$$\mathbb{P}\left\{\max_{0\leq i\leq k-1}\max_{2^{k-1}<n\leq 2^k} n|M_n(x+u_i) - M_n(x)| > Ct_k\right\}$$

(39)
$$\leq \sum_{i=0}^{k-1}\mathbb{P}\left\{\max_{2^{k-1}<n\leq 2^k} n|M_n(x+u_i) - M_n(x)| > Ct_k\right\}$$

$$\leq 2k\exp\left[\frac{-C^2\log k}{2c+1}\right]$$

and similarly

$$\mathbb{P}\left\{\max_{0\leq m\leq 2^k-1}\max_{2^{k-1}<n\leq 2^k} n|M_n(x+v_m) - M_n(x)| > Ct_k\right\}$$

(40)
$$\leq \sum_{m=0}^{2^k-1} 2\exp[-C^2 t_k^2/(2Ct_k + 2\times 2^k cv_m)]$$

$$\leq \sum_{m=0}^{2^k-1} 2\exp[-C^2 t_k^2/(2Ct_k + 2\times 2^k cv_{2^k})]$$

$$\leq 2^{k+1}\exp\left[\frac{-C^2 k\log k}{2c+1}\right].$$

For any $v$, $v_m < v \leq v_{m+1}$, observe that $0 \leq F_n^*(x+v_{m+1}) - F_n^*(x+v_m) \leq cb_{2^k}/(k2^k)$,

$$M_n(x+v) - M_n(x) \leq M_n(x+v_{m+1}) - M_n(x) + cb_{2^k}/(k2^k)$$

and similarly, $M_n(x+v) - M_n(x) \geq M_n(x+v_m) - M_n(x) - cb_{2^k}/(k2^k)$. So (40) yields

$$\mathbb{P}\left\{\sup_{0\leq v\leq \alpha_k}\max_{2^{k-1}<n\leq 2^k} n|M_n(x+v) - M_n(x)| > (C+1)t_k\right\}$$

(41)
$$\leq 2^{k+1}\exp\left[\frac{-C^2 k\log k}{2c+1}\right].$$

Since (41) holds for all $x \in \mathbb{R}$, by the triangle inequality, (41) together with (39) implies

$$\mathbb{P}\left\{\max_{0 \leq u \leq b_{2^k}} \max_{2^{k-1} < n \leq 2^k} n|M_n(x+u) - M_n(x)| > (2C+1)t_k\right\}$$

$$\leq \mathbb{P}\left\{\max_{0 \leq i \leq k-1} \max_{2^{k-1} < n \leq 2^k} n|M_n(x+u_i) - M_n(x)| > Ct_k\right\}$$

$$+ \sum_{i=0}^{k-1} \mathbb{P}\left\{\sup_{0 \leq v \leq \alpha_k} \max_{2^{k-1} < n \leq 2^k} n|M_n(x+v+u_i) - M_n(x+u_i)|\right.$$

$$\left. > (C+1)t_k\right\}$$

$$\leq 2k \exp\left[\frac{-C^2 \log k}{2c+1}\right] + k \times 2^{k+1} \exp\left[\frac{-C^2 k \log k}{2c+1}\right].$$

Therefore (37) follows by letting $C_\tau = 1 + 2(\tau+1)^{1/2}(2c+1)^{1/2}$. $\square$

LEMMA 6. *Assume that the conditions of Lemma 5 are satisfied and in addition assume that there is a $\rho \geq 1$ such that for all sufficiently large $n$ we have that*

(42) $$\frac{b_{2n}}{\rho} \leq \min_{n \leq j \leq 2n} b_j \leq \max_{n \leq j \leq 2n} b_j \leq \rho b_{2n}.$$

*Then for each fixed $x \in \mathbb{R}$,*

(43) $$\sup_{|u| \leq b_n} |M_n(x+u) - M_n(x)| = O_{\text{a.s.}}\left[\frac{\sqrt{b_n \log \log n}}{\sqrt{n}}\right].$$

PROOF. Observe that due to (42), for all sufficiently large $n$ we have

$$\max_{2^{k-1} < n \leq 2^k} \frac{\sqrt{n} \sup_{|u| \leq b_n} |M_n(x+u) - M_n(x)|}{\sqrt{b_n \log \log n}}$$

$$\leq \sup_{|u| \leq \rho b_{2^k}} \max_{2^{k-1} < n \leq 2^k} \frac{n|M_n(x+u) - M_n(x)|}{\sqrt{nb_n \log \log n}}$$

$$\leq \sup_{|u| \leq \rho b_{2^k}} \max_{2^{k-1} < n \leq 2^k} \frac{n|M_n(x+u) - M_n(x)|}{\sqrt{2^{k-1}\rho^{-1}b_{2^k} \log \log 2^{k-1}}}.$$

Hence (43) follows from Lemma 5 via the Borel–Cantelli lemma. $\square$

LEMMA 7. *Assume (4) and that $\mathbb{E}(|X_1|^\alpha) < \infty$ for some $\alpha > 0$. Then for all $\tau > 1$ there exists $C_\tau > 0$ such that*

(44) $$\mathbb{P}\left[\sup_x |M_n(x)| > \frac{C_\tau (\log n)^{1/2}}{n^{1/2}}\right] = O(n^{-\tau})$$

and

$$(45) \quad \mathbb{P}\left[\sup_{|x-y|\leq b_n} |M_n(y) - M_n(x)| > C_\tau \frac{b_n^{1/2}(\log n)^{1/2}}{n^{1/2}}\right] = O(n^{-\tau}),$$

where $(b_n)_{n\geq 1}$ is a positive, bounded sequence of real numbers such that $\log n = o(nb_n)$.

PROOF. We only prove (45) since (44) can be similarly proved. Let $c = \sup_x f_\varepsilon(x) < \infty$, $v_n = \sqrt{nb_n \log n}$, $t_n = v_n/n$, $J = n^{(\tau+4)/\alpha}$ and $Y_i(x) = \mathbf{1}_{X_i \leq x} - \mathbb{E}(\mathbf{1}_{X_i \leq x}|\mathcal{F}_{i-1})$. Then

$$I_n := \mathbb{P}\left[\sup_{|x-y|\leq b_n, x\leq -J} |M_n(y) - M_n(x)| > Ct_n\right]$$

$$\leq n\mathbb{P}\left[\sup_{|x-y|\leq b_n, x\leq -J} |Y_1(y) - Y_1(x)| > Ct_n\right]$$

$$\leq n(Ct_n)^{-1}\mathbb{E}\left[\sup_{|x-y|\leq b_n, x\leq -J} |Y_1(y) - Y_1(x)|\right]$$

$$= O(nt_n^{-1})\mathbb{E}\left[\sup_{x\leq -J+b_n} |Y_1(x)|\right]$$

$$= O(n^2 v_n^{-1})(J - b_n)^{-\alpha}\mathbb{E}(|X_1|^\alpha) = O(n^{-1-\tau}),$$

where Markov's inequality is used in the second inequality. Similarly,

$$III_n := \mathbb{P}\left[\sup_{|x-y|\leq b_n, x\geq J} |M_n(y) - M_n(x)| > Ct_n\right] = O(n^{-1-\tau}).$$

Let $x_i = ib_n/n$, $i = -N-1, \ldots, N+1$, where $N = \lfloor Jn/b_n \rfloor$, and

$$II_n := \mathbb{P}\left[\sup_{|x-y|\leq b_n, -J<x<J} |M_n(y) - M_n(x)| > Ct_n\right].$$

Again by Freedman's inequality, for $|x-y|\leq b_n$ and sufficiently large $n$,

$$\mathbb{P}[n|M_n(y) - M_n(x)| > Cv_n] \leq 2\exp[-C^2 v_n^2/(2Cv_n + 2ncb_n)] \leq 2n^{-C^2/(2c+1)}.$$

Thus

$$\mathbb{P}\left[\max_{i,j=-N-1,\ldots,N+1: |x_i-x_j|\leq b_n} n|M_n(x_i) - M_n(x_j)| > Cv_n\right] = O(N^2)n^{-C^2/(2c+1)}.$$

For any $x, y$ with $|x-y| \leq b_n$, $|x| \leq J$ and $|y| \leq J$, choose $i$ and $j$ such that $x_i \leq x < x_{i+1}$ and $x_j \leq y < x_{j+1}$. Then

$$n[M_n(x_j) - M_n(x_{i+1})] - 2cb_n \leq n[M_n(y) - M_n(x)]$$
$$\leq n[M_n(x_{j+1}) - M_n(x_i)] + 2cb_n.$$

Therefore (45) follows by choosing $C_\tau^2 = (2c+1)[(8+2\tau)/\alpha + \tau + 5]$, given that

$$\mathbb{P}\left\{\sup_{|x-y|\leq b_n} |M_n(y) - M_n(x)| > (C_\tau + 1)t_n\right\} \leq I_n + II_n + III_n$$

by the triangle inequality. $\square$

REMARK 8. In Lemmas 5–7 it is not required that $b_n \to 0$. We shall use this fact to derive (54), which is a key step in proving Theorem 2.

REMARK 9. It is worth noting that Lemmas 5–7 also apply to LRD processes. In Section 7 we will use them to prove the Bahadur representation for LRD processes. For i.i.d. random variables the increments of the empirical and quantile processes are discussed in great detail in [5].

### 6.3. The differentiable part $N_n$.

LEMMA 8. Let $b_n \to 0$. Assume (4) and $\mathbb{E}(|\varepsilon_k|^\alpha) < \infty$ for some $\alpha > 0$. Further assume (5) if $q > 2$ or (8) if $q = 2$. Then

$$(46) \qquad \sup_{|t|\leq b_n} |N_n(x+t) - N_n(x)| = \frac{\ell_q(n)}{\sqrt{n}} O_{\text{a.s.}}(b_n) + O_{\text{a.s.}}(b_n^2).$$

PROOF. Let $c_0 = \sup_x |f'_\varepsilon(x)|$ and recall $f_n^*(x) = dF_n^*(x)/dx$. Clearly $|f'(x)| \leq c_0$ since $f'(x) = \mathbb{E}[f'_\varepsilon(x - X_{i,i-1})]$. Using Taylor's expansion, we get

$$\sup_{|t|\leq b_n} |N_n(x+t) - N_n(x) - t[f_n^*(x) - f(x)]| \leq \frac{b_n^2}{2} \sup_x |d[f_n^*(x) - f(x)]/dx|$$

$$\leq b_n^2 c_0.$$

Let $S_n(x) = n[f_n^*(x) - f(x)]$. If $q > 2$, by (50) of Lemma 10 there exists $C < \infty$ such that $\limsup_{n\to\infty} |S_n(x)|/\sqrt{2n \log\log n} \leq C < \infty$ almost surely. Hence (46) follows. The case that $q = 2$ similarly follows from (ii) of Lemma 10. $\square$

LEMMA 9. Assume (4), (10) and (11). Then for any $-\infty < l < u < \infty$, we have

$$(47) \qquad \mathbb{E}\left[\max_{l\leq x\leq u} (|N_n(x)|^q + |N'_n(x)|^q)\right] = O(n^{-q/2})$$

and

$$(48) \qquad \sup_{x\in[l,u]} [|N_n(x)| + |N'_n(x)|] = o_{\text{a.s.}}[\iota_q(n)/\sqrt{n}].$$

PROOF. We only consider $q > 2$ since the case $q = 2$ can be similarly handled. By Lemma 3 and (34) of Proposition 1, (11) entails $\max_{l \leq x \leq u} \|N_n(x)\|_q = O(1/\sqrt{n})$ and $\max_{l \leq x \leq u} \|N'_n(x)\|_q = O(1/\sqrt{n})$. Since $N_n(x) = N_n(l) + \int_l^x N'_n(t)\,dt$,

$$\mathbb{E}\left[\max_{l \leq x \leq u} |N_n(x)|^q\right] = O\{\mathbb{E}[|N_n(l)|^q]\} + O\left\{\mathbb{E}\left[\int_l^u |N'_n(x)|\,dx\right]^q\right\}$$

$$= O(n^{-q/2}) + O\left\{(u-l)\mathbb{E}\int_l^u |N'_n(x)|^q\,dx\right\}$$

$$= O(n^{-q/2}).$$

Similarly,

$$\mathbb{E}\left[\max_{l \leq x \leq u} |N'_n(x)|^q\right] = O(n^{-q/2}).$$

Then (47) follows. Let $G_n(x) = n[F_n^*(x) - F(x)]$. By Lemma 4, (47) implies that

$$\sum_{k=4}^{\infty} \frac{\mathbb{E}[\max_{n \leq 2^k} \max_{l \leq x \leq u} |G_n(x)|^q]}{2^{qk/2} \iota_q^q(2^k)} = \sum_{k=4}^{\infty} \frac{O[\sum_{j=0}^k 2^{(k-j)/q} 2^{j/2}]^q}{2^{qk/2} \iota_q^q(2^k)}$$

$$= \sum_{k=4}^{\infty} O\{[\iota_q(2^k)]^{-q}\}$$

$$= \sum_{k=4}^{\infty} \frac{O(1)}{k(\log k)^2} < \infty.$$

Then by the Borel–Cantelli lemma, $\max_{l \leq x \leq u} |G_n(x)| = o_{\text{a.s.}}[\iota_q(n)\sqrt{n}]$, which in conjunction with the similar claim $\max_{l \leq x \leq u} |G'_n(x)| = o_{\text{a.s.}}[\iota_q(n)\sqrt{n}]$ entails (48). □

6.4. *Limit theorems for* $F_n - F$.

LEMMA 10. (i) *Assume* (4) *and* (5) *for some* $q > 2$. *Then for every* $x$ *there exist* $0 \leq \sigma_1, \sigma_2 < \infty$ *such that*

(49) $$\limsup_{n \to \infty} \pm \frac{\sqrt{n}[F_n(x) - F(x)]}{\sqrt{2 \log \log n}} = \sigma_1$$

*and*

(50) $$\limsup_{n \to \infty} \pm \frac{\sqrt{n}[f_n^*(x) - f(x)]}{\sqrt{2 \log \log n}} = \sigma_2$$

*almost surely for either choice of sign.*

(ii) *Assume* (4) *and* (8). *Then for every* $x$

(51) $$|F_n(x) - F(x)| + |f_n^*(x) - f(x)| = o_{\text{a.s.}}[\ell_2(n)/\sqrt{n}].$$

PROOF. (i) It is a direct consequence of (35) of Proposition 1 and Lemma 3.

(ii) Let $R_n(x) = n[F_n(x) - F(x)]$. By (8) and (34) of Proposition 1, $\|R_n(x)\| = O(\sqrt{n})$. Then by Lemma 4,

$$\sum_{k=4}^{\infty} \frac{\mathbb{E}[\max_{n \leq 2^k} |R_n(x)|^2]}{2^k \ell_2^2(2^k)} = \sum_{k=4}^{\infty} \frac{O[\sum_{j=0}^{k} 2^{(k-j)/2} 2^{j/2}]^2}{2^k k^3 (\log k)^2}$$

$$= \sum_{k=4}^{\infty} \frac{O(k^2 2^k)}{2^k k^3 (\log k)^2} < \infty,$$

which entails $|R_n(x)| = o_{\text{a.s.}}[\ell_2(n)\sqrt{n}]$ by the Borel–Cantelli lemma. That $|f_n^*(x) - f(x)| = o_{\text{a.s.}}[\ell_2(n)/\sqrt{n}]$ similarly follows. □

LEMMA 11. *Let $q = 2$ and assume that the conditions of Theorem 2 are satisfied. Then $\{\sqrt{n}[F_n(x) - F(x)], \xi_{p_0} \leq x \leq \xi_{p_1}\} \Rightarrow \{W(x), \xi_{p_0} \leq x \leq \xi_{p_1}\}$ for some centered Gaussian process $W$ in the Skorohod space $D[\xi_{p_0}, \xi_{p_1}]$.*

PROOF. It suffices to verify the finite-dimensional convergence and the tightness [2]. By Lemma 3 $\|\mathcal{P}_0 \mathbf{1}_{X_n \leq x}\| = O[|a_n|^{\min(\alpha/2, 1)}]$, which is summable in view of (11) since $q = 2$. Then by the Cramér–Wold device, the finite-dimensional convergence easily follows from Lemma 3 in [31].

Write $l = \xi_{p_0}$ and $u = \xi_{p_1}$. Recall $F_n(x) - F(x) = M_n(x) + N_n(x)$. To show the tightness of $\{\sqrt{n}[F_n(x) - F(x)], l \leq x \leq u\}$, it suffices to show that both $\{\sqrt{n} N_n(x), l \leq x \leq u\}$ and $\{\sqrt{n} M_n(x), l \leq x \leq u\}$ are tight. The former easily follows from

$$\mathbb{E}\left[\sup_{|x-y| \leq \delta, l \leq x, y \leq u} n|N_n(x) - N_n(y)|^2\right] \leq \delta^2 n \mathbb{E}\left[\sup_{l \leq \theta \leq u} |f_n^*(\theta) - f(\theta)|^2\right] \leq C\delta^2$$

in view of (47) of Lemma 9 with $q = 2$. For the latter, let $d_i = \mathbf{1}_{x < X_i \leq y} - \mathbb{E}(\mathbf{1}_{x < X_i \leq y} | \mathcal{F}_{i-1})$, $l \leq x < y \leq u$. Then by (4) $\mathbb{E}(d_i^2 | \mathcal{F}_{i-1}) \leq C(y - x)$. Here $C$ denotes a constant which does not depend on $n$, $x$ and $y$ and it may vary from line to line. By Burkholder's inequality [3],

$$\mathbb{E}[n^2 |M_n(x) - M_n(y)|^4] \leq \frac{C}{n^2} \left\|\sum_{i=1}^{n} d_i^2\right\|^2$$

$$\leq \frac{C}{n^2} \left\|\sum_{i=1}^{n} (d_i^2 - \mathbb{E}(d_i^2 | \mathcal{F}_{i-1}))\right\|^2 + \frac{C}{n^2} \|\mathbb{E}(d_i^2 | \mathcal{F}_{i-1})\|^2$$

$$\leq \frac{C}{n} \|d_1^2 - \mathbb{E}(d_1^2 | \mathcal{F}_0)\|^2 + C(y - x)^2$$

$$\leq \frac{C}{n}(y - x) + C(y - x)^2.$$

See inequality (48) in [31] for a similar claim. Therefore, by the argument of Theorem 22.1 in [2], pages 197–199, the process $\{\sqrt{n}M_n(x), l \leq x \leq u\}$ is tight. □

REMARK 10. Under conditions of the type given in (8), Wu [32] obtained a central limit theorem for $S_n(K)/\sqrt{n}$, where $S_n(K) = \sum_{i=1}^{n}[K(X_i) - \mathbb{E}K(X_i)]$, $K$ is a measurable function and $\varepsilon_i$ may have infinite variance.

LEMMA 12. Let $X_k = \sum_{i=0}^{\infty} a_i \varepsilon_{k-i}$ and assume (4) and $\mathbb{E}(|\varepsilon_k|^\alpha) < \infty$ for some $\alpha > 0$. Further assume (42) and $\log^3 n = o(nb_n)$.

(i) If (5) holds with $q > 2$, then for every fixed $x$,

$$\sup_{|u| \leq b_n} |F_n(x+u) - F(x+u) - [F_n(x) - F(x)]|$$

(52)
$$= \frac{O_{\text{a.s.}}(\sqrt{b_n \log \log n})}{\sqrt{n}} + \frac{O_{\text{a.s.}}[b_n \ell_q(n)]}{\sqrt{n}} + O_{\text{a.s.}}(b_n^2).$$

(ii) If (8) holds, then we have (52) with $q = 2$.

REMARK 11. The second term $O_{\text{a.s.}}[b_n \ell_q(n)]/\sqrt{n}$ in the bound of (52) is needed only when $q = 2$.

Lemma 12 follows from Lemmas 6 and 8 and it provides a local fluctuation rate of empirical processes for linear processes. The last two terms of (52) are due to the presence of dependence, in the sense that they disappear if $X_k$ are i.i.d. Actually, if $X_i$ are i.i.d., then $F_n^* \equiv F$ and hence $N_n \equiv 0$.

LEMMA 13. Assume (4), (10) and (11). Then under the conditions of Lemma 7, we have for any $-\infty < l < u < \infty$ that

$$\sup_{|x-y| \leq b_n, \ x,y \in [l,u]} |[F_n(x) - F(x)] - [F_n(y) - F(y)]|$$

(53)
$$= O_{\text{a.s.}}\left[\frac{\sqrt{b_n \log n}}{\sqrt{n}} + \frac{b_n \iota_q(n)}{\sqrt{n}}\right].$$

PROOF. By Lemma 7 it suffices to show that

$$\sup_{|x-y| \leq b_n, \ x,y \in [l,u]} |[F_n^*(x) - F(x)] - [F_n^*(y) - F(y)]| \leq b_n \sup_{\theta \in [l,u]} |f_n^*(\theta) - f(\theta)|$$

$$= b_n o_{\text{a.s.}}[\iota_q(n)/\sqrt{n}],$$

which is an easy consequence of Lemma 9. □

6.5. *Proofs.*

PROOF OF THEOREM 1. We only consider $q > 2$ since the case $q = 2$ follows along similar lines.

(i) Let $b_n = \delta_{n,q}$. Then (42) holds. By Lemma 12 there exists a constant $C_1 < \infty$ such that

$$|[F_n(\xi_p + b_n) - F(\xi_p + b_n)] - [F_n(\xi_p) - F(\xi_p)]| \leq C_1\sqrt{(b_n \log\log n)/n}$$

almost surely. Observe that $F(\xi_p + b_n) = F(\xi_p) + b_n f(\xi_p) + O(b_n^2)$ in view of (4) via Taylor's expansion. By (i) of Lemma 10, there exists a constant $C_2 < \infty$ such that $n|F_n(\xi_p) - F(\xi_p)| \leq C_2 \sqrt{n} \ell_q(n)$ almost surely. Choose $C > 0$ such that $C - C_2 - C_1\sqrt{C/f(\xi_p)} \geq 1$, namely, $C \geq [C_1/\sqrt{f(\xi_p)} + \sqrt{C_1^2/f(\xi_p) + 4(1 + C_2)}]^2/4$. Then for $b_n = C\ell_q(n)/[f(\xi_p)\sqrt{n}]$, $F_n(\xi_p + b_n) > p$ holds almost surely. The other statement that $p > F_n(\xi_p - b_n)$ almost surely similarly follows. Let $\Delta_n = \xi_{n,p} - \xi_p$. Since $F_n$ is nondecreasing, by (6) $|\Delta_n| \leq b_n$ almost surely.

(ii) The argument for Theorem 4 can be applied here. Applying Lemma 12 with $x = \xi_p$, we have

$$|F_n(\xi_{n,p}) - F(\xi_p + \Delta_n) - [F_n(\xi_p) - F(\xi_p)]| = O_{\text{a.s.}}[\sqrt{(b_n \log\log n)/n}\,].$$

Notice that $|F_n(\xi_{n,p}) - p| \leq 1/n$ and, by Taylor's expansion $F(\xi_p + \Delta_n) = p + \Delta_n f(\xi_p) + O(\Delta_n^2)$ since $\sup_x |f'(x)| < \infty$. Then

$$\Delta_n f(\xi_p) = p - F_n(\xi_p) + O_{\text{a.s.}}[\sqrt{(b_n \log\log n)/n}\,]$$

and it entails (7). □

PROOF OF THEOREM 2. Let $l = \xi_{p_0}$ and $u = \xi_{p_1}$. By Lemma 6 and (48) of Lemma 9, we have

$$\sup_{x \in [l,u]} |F_n(x) - F(x)| \leq \sup_{x \in [l,u]} |F_n(x) - F_n^*(x)| + \sup_{x \in [l,u]} |F_n^*(x) - F(x)|$$

(54)
$$= O_{\text{a.s.}}\left[\frac{\sqrt{\log\log n}}{\sqrt{n}}\right] + o_{\text{a.s.}}\left[\frac{\iota_q(n)}{\sqrt{n}}\right] = o_{\text{a.s.}}\left[\frac{\iota_q(n)}{\sqrt{n}}\right].$$

Let $b_n = \iota_q(n)/\sqrt{n}$. (i) By Lemma 13,

$$\inf_{l \leq x \leq u}[F_n(x + b_n) - F(x)]$$

$$\geq \inf_{l \leq x \leq u}[F(x + b_n) - F(x)]$$

$$- \sup_{l \leq x \leq u} |F_n(x) - F(x)|$$

$$- \sup_{|x-y|\leq b_n,\ l\leq x,y\leq u} |[F_n(x) - F(x)] - [F_n(y) - F(y)]|$$

$$\geq b_n \inf_{p_0 \leq p u} f(\xi_p) + O(b_n^2) + o_{\text{a.s.}}(b_n)$$

$$+ O_{\text{a.s.}}[\sqrt{b_n(\log n)/n} + b_n \iota_q(n)/\sqrt{n}].$$

Hence $\inf\{F_n(x + b_n) - F(x): l \leq x \leq u\} > 0$ almost surely, which implies (i) together with a similar claim that $\sup\{F_n(x - b_n) - F(x): l \leq x \leq u\} < 0$ almost surely. The representation (12) then follows from Lemma 13 by using the same argument as in the proof of (ii) of Theorem 1. □

**7. Proof and the sharpness of Theorem 3.** In the study of LRD processes, the asymptotic expansion of empirical processes plays an important role [17, 31]. Let $U_{n,r} = \sum_{0 \leq j_1 < \cdots < j_r} \prod_{s=1}^r a_{j_s} \varepsilon_{n-j_s}$, $U_{n,0} = 1$. For a nonnegative integer $\rho$, similarly to (4) in [31] let

$$S_n(y; \rho) = \sum_{i=1}^n \left[ \mathbf{1}(X_i \leq y) - \sum_{r=0}^{\rho} (-1)^r F^{(r)}(y) U_{i,r} \right];$$

see also [17]. The quantity $S_n(y; \rho)$ can be viewed as the remainder of the $\rho$th-order expansion of $F_n(y)$. In our derivation of Bahadur's representation for LRD processes, we only deal with $\rho = 1$ and do not pursue the higher-order case $\rho \geq 2$ since it involves some really cumbersome manipulations.

As in [31], let $\theta_n = |a_{n-1}|[|a_{n-1}| + (\sum_{i=n-1}^\infty a_i^2)^{1/2} + (\sum_{i=n-1}^\infty a_i^4)^{\rho/2}]$, $\Theta_n = \sum_{i=1}^n \theta_i$, $\Xi_n = n\Theta_n^2 + \sum_{i=1}^\infty (\Theta_{n+i} - \Theta_i)^2$. Since $\rho = 1$, $\theta_n = O[|a_{n-1}| \times (\sum_{i=n-1}^\infty a_i^2)^{1/2}]$. Recall that $\Psi_n = \sqrt{n} \sum_{k=1}^n k^{1/2 - 2\beta} L^2(k)$, $A_n(\beta) = \Psi_n^2 (\log n) \times (\log \log n)^2$ if $\beta < 3/4$ and $A_n(\beta) = \Psi_n^2 (\log n)^3 (\log \log n)^2$ if $\beta \geq 3/4$. Let $H_n(y) = n[F_n^*(y) - F(y) + f(y)\bar{X}_n]$ and $h_n(y) = dH_n(y)/dy$.

LEMMA 14. *Assume $\mathbb{E}(\varepsilon_i^4) < \infty$ and*

$$(55) \qquad \sup_x |f_\varepsilon(x)| + \sup_x |f'_\varepsilon(x)| + \int_{\mathbb{R}} |f'_\varepsilon(u)|^2\, du < \infty.$$

*Then*

$$(56) \qquad \left\| \sup_y |H_n(y)| \right\| + \left\| \sup_y |h_n(y)| \right\| = O(\Psi_n).$$

PROOF. Let $I = \int_R |f'_\varepsilon(u)|^2\, du$ and $K_\theta(x) = [f_\varepsilon(\theta - x) - f_\varepsilon(\theta)]/\sqrt{I}$. Then $k_\theta(x) = \partial K_\theta(x)/\partial x = -f'_\varepsilon(\theta - x)/\sqrt{I}$ satisfies $\int_{\mathbb{R}} k_\theta^2(x)\, dx = 1$. Hence for all $\theta$,

$$K_\theta \in \mathcal{K}(0) := \left\{ K(x) = \int_0^x g(t)\, dt : \int_{\mathbb{R}} g^2(t) \leq 1 \right\};$$

see [31] for the definition of the class $\mathcal{K}$. By Theorem 1 in [31], for

$$S_n(K_\theta, 1) = \frac{1}{\sqrt{I}} \sum_{i=1}^n [f_\varepsilon(\theta - X_{i,i-1}) - f(\theta) + f'(\theta) X_{i,i-1}]$$

we have

$$\mathbb{E}\left[\sup_{\theta \in \mathbb{R}} S_n^2(K_\theta, 1)\right] = \mathcal{O}(\Xi_n).$$

Notice that $S_n(K_\theta, 1)\sqrt{I} - h_n(\theta) = -f'(\theta) \sum_{i=1}^n \varepsilon_i$. Then $\|\sup_y |h_n(y)|\| = O(\Xi_n^{1/2})$ since $\sup_\theta |f'(\theta)| < \infty$ and $\|\sum_{i=1}^n X_{i,i-1} - n\bar{X}_n\| = O(\sqrt{n})$. By Karamata's theorem, it is easily seen that $\Xi_n = O(\Psi_n^2)$ (cf. Lemma 5 in [31]). Similarly, $\|\sup_y |H_n(y)|\| = O(\Xi_n^{1/2})$ holds under the condition $\int_\mathbb{R} f_\varepsilon^2(u)\,du < \infty$. The last inequality trivially holds since $\sup_u f_\varepsilon(u) < \infty$. $\square$

LEMMA 15.  *Assume* $\mathbb{E}(\varepsilon_i^4) < \infty$ *and* (55).

(i) *Let* $(\delta_n)_{n \geq 1}$ *be a positive, bounded sequence such that* $\log n = o(n\delta_n)$. *Then*

(57) $$\sup_{|x-y| \leq \delta_n} |S_n(y; 1) - S_n(x; 1)| = O_{\text{a.s.}}[\sqrt{n\delta_n \log n} + \delta_n A_n^{1/2}(\beta)].$$

(ii) *For any* $-\infty < l < u < \infty$, $\sup_{l \leq y \leq u} |S_n(y; 1)| = o_{\text{a.s.}}[A_n^{1/2}(\beta)]$.

PROOF. (i) By Lemma 7, since $\mathbb{E}(X_1^2) < \infty$, $\sqrt{n}\sup_{|x-y| \leq \delta_n} |M_n(y) - M_n(x)| = O_{\text{a.s.}}(\sqrt{\delta_n \log n})$. To show (57), notice that $\sup_{|x-y| \leq \delta_n} |H_n(y) - H_n(x)| \leq \delta_n \sup_\theta |h_n(\theta)|$; it suffices to verify that $\sup_\theta |h_n(\theta)| = o_{\text{a.s.}}[A_n^{1/2}(\beta)]$ in view of

(58) $$S_n(y; 1) - S_n(x; 1) = n[M_n(y) - M_n(x)] + [H_n(y) - H_n(x)].$$

By Karamata's theorem, $\sum_{j=0}^d 2^{(d-j)/2} \Psi_{2^j} = O(\Psi_{2^d})$ if $\beta < 3/4$ and $\sum_{j=0}^d 2^{(d-j)/2} \Psi_{2^j} = O(d\Psi_{2^d})$ if $\beta \geq 3/4$. So it follows from Lemma 14 that

$$\sum_{j=0}^d 2^{(d-j)/2} \left\|\sup_y |h_{2^j}(y)|\right\| = \sum_{j=0}^d 2^{(d-j)/2} O(\Psi_{2^j}) = \frac{O[A_{2^d}^{1/2}(\beta)]}{d^{1/2} \log d},$$

which in conjunction with Lemma 4 implies

$$\sum_{d=3}^\infty \frac{1}{A_{2^d}(\beta)} \mathbb{E}\left[\max_{j \leq 2^d} \sup_y |h_j(y)|^2\right] = \sum_{d=3}^\infty O(d^{-1} \log^{-2} d) < \infty.$$

Hence $\sup_y |h_n(y)| = o_{\text{a.s.}}[\sqrt{A_n(\beta)}]$ via the Borel–Cantelli lemma.

(ii) Notice that $S_n(y;1) = nM_n(y) + H_n(y)$. By Lemma 7,
$$\sqrt{n} \sup_{l \leq y \leq u} |M_n(y)| = O_{\text{a.s.}}(\sqrt{\log n}).$$

Using the argument in (i), (56) implies $\sup_y |H_n(y)| = o_{\text{a.s.}}[\sqrt{A_n(\beta)}]$. Hence (ii) follows in view of $\sqrt{n} = O(\Psi_n)$ and $\sqrt{n \log n} = o[\sqrt{A_n(\beta)}]$. □

LEMMA 16. *Assume $\mathbb{E}(\varepsilon_i^4) < \infty$ and (55). Let $B_n = \sigma_{n,1}(\log n)^{1/2}(\log \log n)$, $b_n = B_n/n$ and $\Delta_{n,p} = \xi_{n,p} - \xi_p$. Then (i) $\bar{X}_n = o_{\text{a.s.}}(b_n)$ and (ii) if, in addition, $\inf_{p_0 \leq p \leq p_1} f(\xi_p) > 0$ for some $0 < p_0 < p_1 < 1$, we have*

(59) $$\sup_{p_0 \leq p \leq p_1} |\Delta_{n,p}| = o_{\text{a.s.}}(b_n)$$

*and*

(60) $$\sup_{p_0 \leq p \leq p_1} |\Delta_{n,p} - \bar{X}_n| = o_{\text{a.s.}}(b_n^2) + o_{\text{a.s.}}[n^{-1} A_n^{1/2}(\beta)].$$

PROOF. (i) Let $S_n = \sum_{i=1}^n X_i$. Since $\sigma_{n,1} = \|S_n\| \sim Cn^{3/2-\beta}L(n)$, by Lemma 4
$$B_{2^d}^{-2} \left\|\max_{i \leq 2^d} |S_i|\right\|^2 \leq B_{2^d}^{-2} \left[\sum_{r=0}^d 2^{(d-r)/2} \sigma_{2^r,1}\right]^2 = O(d^{-1} \log^{-2} d).$$

Again by the Borel–Cantelli lemma $\bar{X}_n = o_{\text{a.s.}}(b_n)$.

(ii) Similarly as in the proof of Theorem 1, it suffices to show that, due to the monotonicity of $F_n$, $\inf_{p_0 \leq p \leq p_1}[F_n(\xi_p + b_n) - p] > 0$ holds almost surely since the other inequality $\sup_{p_0 \leq p \leq p_1}[F_n(\xi_p - b_n) - p] < 0$ can be similarly derived. By Lemma 15

$$\inf_{p_0 \leq p \leq p_1}[F_n(\xi_p + b_n) - p]$$
$$\geq \inf_{p_0 \leq p \leq p_1}[F(\xi_p + b_n) - f(\xi_p + b_n)\bar{X}_n - p + S_n(\xi_p;1)/n]$$
$$\quad - \sup_{|x-y| \leq b_n} |S_n(y;1) - S_n(x;1)|/n =: I_n + II_n.$$

Since $\sup_x |f'_\varepsilon(x)| < \infty$, by Taylor's expansion $\sup_p |F(\xi_p + b_n) - p - b_n f(\xi_p)| = O(b_n^2)$ and $\sup_p |f(\xi_p + b_n) - f(\xi_p)| = O(b_n)$. Let $l = \xi_{p_0}$ and $u = \xi_{p_1}$. By (ii) of Lemma 15, $\sup_{l \leq x \leq u} |S_n(x;1)| = o_{\text{a.s.}}[A_n^{1/2}(\beta)]$. By (i) $\bar{X}_n = o_{\text{a.s.}}(b_n)$. Therefore

$$I_n = \inf_{p_0 \leq p \leq p_1} f(\xi_p)(b_n - \bar{X}_n) + O(b_n^2 + b_n|\bar{X}_n|) + o_{\text{a.s.}}[A_n^{1/2}(\beta)/n]$$
$$\geq \tfrac{1}{2} \inf_{p_0 \leq p \leq p_1} f(\xi_p) b_n$$

almost surely. By (57) of Lemma 15, $II_n = o_{\text{a.s.}}(b_n)$ and hence (59) holds. Relation (60) follows by letting $y = \xi_{n,p} = \xi_p + \Delta_{n,p}$ in (ii) of Lemma 15 in view of

$$\sup_{p_0 \leq p \leq p_1} |F(\xi_p + \Delta_{n,p}) - p - f(\xi_p)\Delta_{n,p}|$$
$$\leq \frac{\sup_x |f'(x)|}{2} \sup_{p_0 \leq p \leq p_1} \Delta_{n,p}^2 = o_{\text{a.s.}}(b_n^2)$$

and $\sup_{p_0 \leq p \leq p_1} |f(\xi_p + \Delta_{n,p}) - f(\xi_p)| = o_{\text{a.s.}}(b_n)$. □

REMARK 12. Under the stronger condition that $f_\varepsilon$ is four times differentiable with bounded, continuous and integrable derivatives, Ho and Hsing [17] obtained

(61) $$\sup_{p_0 \leq p \leq p_1} |\Delta_{n,p} - \bar{X}_n| = o_{\text{a.s.}}(n^{-1-\lambda}\sigma_{n,1})$$

for all $0 < \lambda < \min(1 - \beta, \beta - 1/2)$; see Theorem 5.1 therein. The result (61) is very interesting in the sense that $\Delta_{n,p}$ can be approximated by $\bar{X}_n$, which *does not* depend on $p$. Consequently the asymptotic distribution of the trimmed and Winsorized means easily follows from that of $\bar{X}_n$. After elementary calculations it is easily seen that our bound (60) is slightly sharper.

PROOF OF THEOREM 3. By (59) $\sup_{p_0 \leq p \leq p_1} |\Delta_{n,p}| = o_{\text{a.s.}}(b_n)$. Applying Lemma 15 with $x = \xi_p$ and $y = \xi_{n,p}$, $p_0 \leq p \leq p_1$, we have

(62)
$$n \sup_{p_0 \leq p \leq p_1} |p - F(\xi_p + \Delta_{n,p}) + f(\xi_{n,p})\bar{X}_n - [F_n(\xi_p) - F(\xi_p) + f(\xi_p)\bar{X}_n]|$$
$$= O_{\text{a.s.}}[\sqrt{nb_n \log n} + b_n A_n^{1/2}(\beta)].$$

Since $\sup_x[|f'(x)| + |f''(x)|] < \infty$, by Taylor's expansion

$$\sup_{p_0 \leq p \leq p_1} |F(\xi_p + \Delta_{n,p}) - p - \Delta_{n,p} f(\xi_p) - \Delta_{n,p}^2 f'(\xi_p)/2| = o_{\text{a.s.}}(b_n^3)$$

and

$$\sup_{p_0 \leq p \leq p_1} |f(\xi_p + \Delta_{n,p}) - f(\xi_p) - \Delta_{n,p} f'(\xi_p)| = o_{\text{a.s.}}(b_n^2).$$

After some elementary calculations, (62) implies

$$\sup_{p_0 \leq p \leq p_1} \left| f(\xi_p)\Delta_{n,p} + \frac{f'(\xi_p)}{2}(\Delta_{n,p} - \bar{X}_n)^2 - \frac{1}{2}f'(\xi_p)\bar{X}_n^2 - [p - F_n(\xi_p)] \right|$$
$$= o_{\text{a.s.}}(b_n^3) + n^{-1} O_{\text{a.s.}}[\sqrt{nb_n \log n} + b_n A_n^{1/2}(\beta)].$$

Observe that $\Psi_n = O[\sqrt{n}L^*(n) + n^{2-2\beta}L^2(n)]$ and $A_n^{1/2}(\beta) \leq \Psi_n (\log n)^{3/2} \times (\log \log n) = o(nb_n)$. Thus (15) follows from (60) and

$$\sup_{p_0 \leq p \leq p_1} (\Delta_{n,p} - \bar{X}_n)^2 = o_{\text{a.s.}}[b_n^2 + A_n^{1/2}(\beta)/n]^2$$
$$= o_{\text{a.s.}}[b_n^4 + A_n(\beta)/n^2]$$
$$= o_{\text{a.s.}}[b_n^3 + b_n A_n^{1/2}(\beta)/n]. \qquad \square$$

7.1. *The sharpness of Theorem 3.* It is challenging to obtain a sharp bound for the left-hand side of (15) in Theorem 3. We now comment on the sharpness of Lemma 15, which describes the oscillations of $F_n(x) - F(x) + f(x)\bar{X}_n$. Recall that in the SRD case the sharp oscillation rate of $F_n(x) - F(x)$ at a fixed $x$ in Lemma 12 leads to the Bahadur representation with optimal bound by letting $b_n = c\sqrt{(\log \log n)/n}$ for some $c > 0$. Here we claim that the bound in (57) of Lemma 15, which is a key ingredient for the derivation of (15), is optimal up to a multiplicative slowly varying function.

LEMMA 17. *Assume $\mathbb{E}(\varepsilon_i^4) < \infty$, (14) and $\int_{\mathbb{R}} |f_\varepsilon''(u)|^2\, du < \infty$. Let $\delta_n = n^\gamma L_2(n)$ for some slowly varying function $L_2$, $-1 < \gamma < 0$ and $\sigma_{n,2} = n^{2-2\beta} L^2(n)$.*

(i) *If $4\beta - 3 > \gamma$, then $[S_n(x + \delta_n; 1) - S_n(x; 1)]/\sqrt{n\delta_n} \Rightarrow N[0, f(x)]$.*
(ii) *If $4\beta - 3 < \gamma$, then*

(63)
$$\frac{S_n(x + \delta_n; 1) - S_n(x; 1)}{\sigma_{n,2} \delta_n}$$
$$\Rightarrow f'(x) C_\beta \int_{u_1 < u_2 < 1} \int_0^1 [(v - u_1)_+ (v - u_2)_+]^{-\beta}\, dv\, d\mathbb{B}(u_1)\, d\mathbb{B}(u_2)$$

*for some constant $C_\beta > 0$, where $\mathbb{B}$ is a standard two-sided Brownian motion and $z_+ = \max(z, 0)$. In particular, if $\gamma = 1/2 - \beta$, then* (i) *[resp.* (ii)*] holds if $7/10 < \beta < 1$ [resp. $1/2 < \beta < 7/10$].*

REMARK 13. The limiting distribution in (63) is called the Rosenblatt distribution, a special case of multiple Wiener–Itô integrals [23].

PROOF OF LEMMA 17. Observe that

$$S_n(x + \delta_n; 1) - S_n(x; 1) = n[M_n(x + \delta_n) - M_n(x)] + [H_n(x + \delta_n) - H_n(x)].$$

Write $n[F_n(x+\delta_n) - F_n(x)] = \sum_{i=1}^n K[(x-X_i)/\delta_n]$, where the kernel $K(u) = \mathbf{1}_{-1 \leq u \leq 0}$. By Lemma 2 in [34] $n[M_n(x+\delta_n) - M_n(x)]/\sqrt{n\delta_n} \Rightarrow N[0, \sigma^2(x)]$ with $\sigma^2(x) = f(x) \int_\mathbb{R} K^2(u)\, du = f(x)$. By Lemma 14

$$\|H_n(x+\delta_n) - H_n(x)\| \le \delta_n \left\|\sup_y |h_n(y)|\right\|$$

$$= O(\delta_n \Psi_n) = \delta_n O[\sqrt{n} L^*(n) + n^{2-2\beta} L^2(n)].$$

If $4\beta - 3 > \gamma$, then $\delta_n \Psi_n = o(\sqrt{n\delta_n})$ and (i) follows.

On the other hand, if $4\beta - 3 < \gamma$, then $\beta \in (1/2, 3/4)$ and $h_n(x)/\sigma_{n,2}$ converges to the Rosenblatt distribution in (63); see Lemma 4 in [34] and Corollary 3 in [31]. Under the conditions (14) and $\int_\mathbb{R} |f''_\varepsilon(u)|^2\, du < \infty$, by the argument of Lemma 14, we have $\|\sup_u |h'_n(u)|\| = O(\Psi_n)$. Then $\|H_n(x+\delta_n) - H_n(x) - \delta_n h_n(x)\| \le \frac{1}{2}\delta_n^2 \|\sup_u |h'_n(u)|\| = O(\delta_n^2 \Psi_n)$ and (ii) follows in view of $\sqrt{n\delta_n} + \delta_n^2 \Psi_n = o(\delta_n \sigma_{n,2})$.

If $\gamma = 1/2 - \beta$, then $4\beta - 3 > \gamma$ if and only if $7/10 < \beta$. $\square$

Lemma 17 asserts the dichotomous convergence of $S_n(x+\delta_n; 1) - S_n(x; 1)$ at a fixed point $x$. Notice that $\bar{X}_n/[n^{1/2-\beta}L(n)] \Rightarrow N(0, \sigma^2)$ and, by (60), $(\xi_{n,p} - \xi_p)/[n^{1/2-\beta}L(n)] \Rightarrow N(0, \sigma^2)$ for some $\sigma^2 < \infty$. For $\delta_n = n^{1/2-\beta}L_2(n)$, Lemma 17 shows that, up to a multiplicative slowly varying function, the optimal bound of $[S_n(x+\delta_n; 1) - S_n(x; 1)]/n$ is $n^{\max(-\beta/2-1/4,\ 3/2-3\beta)}$. This bound indicates that (15) or (16) is optimal up to a multiplicative slowly varying function. It also explains why there is a boundary $\beta = 7/10$ in (15) or (16); see the discussion of the three terms in the $O_{\text{a.s.}}$ bound of (15) in Section 3.

**Acknowledgments.** I am grateful to the referees and an Associate Editor for their many helpful comments. I would also like to thank Jan Mielniczuk, George Michailidis and Xiaofeng Shao for valuable suggestions.

Department of Statistics
University of Chicago
5734 S. University Avenue
Chicago, Illinois 60637
USA
e-mail: wbwu@galton.uchicago.edu